\documentclass[11pt]{article}
\usepackage[utf8]{inputenc}
\usepackage[a4paper, margin=1in]{geometry}
\usepackage{graphicx}
\usepackage{amsmath, amssymb}
\usepackage{mhchem}
\usepackage{setspace}
\usepackage{xcolor}
\usepackage{natbib}
\usepackage{hyperref}
\usepackage{titlesec}
\usepackage{float}
\usepackage{tikz}
\usepackage{longtable}  
\usetikzlibrary{positioning, arrows.meta}
\usepackage{appendix}
\titleformat{\section}{\large\bfseries}{\thesection.}{1em}{}
\titleformat{\subsection}{\normalsize\bfseries}{\thesubsection.}{1em}{}

\title{\textbf{A Network-Based Framework to Identify Synergies and Trade-offs among SDG Indicators}}
\author{
  Gaurav Kottari\thanks{Email: gk917@snu.edu.in} \\
  \small Department of Mathematics \\
  \small Shiv Nadar Institution of Eminence\\
  \small Delhi-NCR, India
  \and
  Qazi J. Azhad\thanks{Email: qazi.jamal@snu.edu.in} \\
  \small Department of Mathematics\\
  \small Shiv Nadar Institution of Eminence\\
  \small Delhi-NCR, India
  \and
  Niteesh Sahni\thanks{Corresponding Email: niteesh.sahni@snu.edu.in} \\
  \small Department of Mathematics\\
  \small Shiv Nadar Institution of Eminence\\
  \small Delhi-NCR, India
}
\date{}

\begin{document}

\maketitle

\begin{abstract}
Achieving the United Nations Sustainable Development Goals (SDGs) requires an understanding of the complex interlinkages that exist among their underlying indicators. While most existing research examines these interconnections at the goal level, policy interventions are typically designed and implemented at the indicator level, where synergies and trade-offs most directly emerge. This study addresses this gap by proposing a network-theoretic framework to assess indicator-level interactions in a systematic and data-driven manner. We introduce two complementary measures—the positive strength and negative strength of an indicator, which jointly capture the balance between synergistic and conflicting interactions within a national SDG indicator network. Based on these measures, indicators are classified as synergy- and trade-off-dominated according to their net systemic interaction structure. To move beyond classification, we further examine the structural drivers of synergy dominance using an explanatory regression framework, focusing on the roles of direct positive interactions and indirect network embeddedness. This analysis shows that indicators classified as synergy-dominated are typically characterised by a high concentration of direct synergies and additional support from indirect pathways through the network, allowing positive effects to extend beyond immediate neighbours. The framework is applied to two national case studies, India and Italy, to illustrate how the classification of indicators vary across development contexts. Overall, the proposed methodology provides a transparent and scalable tool for identifying the structural conditions under which indicator-level synergies emerge, thereby supporting a more nuanced understanding of how development actions can generate reinforcing effects across the SDG system.
\end{abstract}


\section{Introduction}
Sustainable Development Goals (SDGs) were adopted by the member states of the United Nations in 2015, which is a universal call to eradicate poverty, protect the environment, and guarantee peace and prosperity for all by 2030~\citep{assembly2015resolution}. It consists of 17 goals, 169 targets, and 234 indicators that address social, economic, and environmental issues. These are very broad and interconnected in nature, as progress on one goal rarely happens in isolation. Gains in areas such as health, education, or clean energy can enable successes in other areas but can also lead to unexpected setbacks or new pressures elsewhere. Hence, identifying these interlinkages is essential to understand the general dynamics of sustainable development.

Since each goal is further divided into targets, and each target breaks down into indicators guiding specific actions and measures, it is important to look into interactions not only among goals but also among their underlying indicators; this is often the level at which the most direct trade-offs and synergies appear. The links between indicators are therefore crucial in finding synergies that accelerate progress or trade-offs that could slow it down~\citep{swain2021modeling}. 

To address this need for understanding interlinkages, existing studies have used two broad approaches: qualitative and quantitative. In qualitative studies, the strength of the interlinkages between targets (or goals) is determined by expert opinions, policy analysis, or literature reviews. Nilsson~\citep{nilsson2016policy} proposed a seven-point interaction scale that helps classify whether one target supports or harms another. Building on this scale, Allen et al.~\citep{allen2019prioritising} and Weitz et al.~\citep{weitz2018towards} explicitly prioritised SDG targets by identifying those that create many positive effects across the goals. Subsequent studies~\citep{xiao2023synergies,fader2018toward} further applied qualitative scoring, refined interaction assessments to map how targets influence one another. Together, these qualitative studies indicate which SDG targets may produce strong positive links and where negative interactions may arise, constituting a well-established body of work that evaluates targets by their systemic roles and cross-goal impacts.

On the other hand, quantitative research uses statistical and network-based methods to analyse SDG interlinkages using empirical indicator data. Several studies~\citep{pradhan2017systematic,de2020synergies,hegre2020synergies,kostetckaia2022sustainable,miao2025priority} have applied correlation analysis, principal component analysis, and clustering techniques to describe synergy and trade-off patterns across goals and regions. Although these contributions have significantly advanced the understanding of SDG interactions, they focus primarily on interactions between goals rather than at the target or indicator level.

Because policy interventions are made at the level of specific indicators, determining which of these are synergistic becomes essential for maximising development gains. By “synergistic”, we mean indicators whose progress tends to produce net positive effects across other indicators, where the positive spillovers from improving these indicators outweigh any negative effects on others. Since interventions on these indicators reinforce progress not only within their immediate domain but also across multiple interconnected areas, identifying such indicators guarantees that available resources are invested in generating the maximum systemic payoffs. In addition, finding indicators that are trade-offs is equally important. If taken carelessly, actions based on these indicators could unintentionally impede or even reverse progress in other areas. By recognising both synergistic and trade-off indicators, policymakers can design strategies that strengthen complementarities while minimising potential conflicts, thereby enhancing the overall efficiency and effectiveness of sustainable development efforts.

Although qualitative studies have advanced target-level prioritisation~\citep{allen2019prioritising,weitz2018towards} and indirectly informed indicator-level assessments, a parallel, data-driven framework for identifying influential indicators within quantitative science remains critically absent. Qualitative approaches rely heavily on expert judgement, literature synthesis, and manual scoring, limiting their ability to be applied systematically across all indicators, countries, and time periods. A quantitative, data-driven approach can overcome these limitations by evaluating indicator interactions using empirical data, enabling reproducible, comparable, and timely identification of influential indicators, particularly at the national level where policy decisions are implemented. 

To the best of current knowledge, only two quantitative studies have attempted to identify significant targets or indicators, and both exhibit critical limitations:
\begin{itemize}
    \item Ranganathan and Swain~\citep{swain2021modeling} proposed a methodology using network theory to prioritise SDG indicators. They constructed regional correlational networks of SDG indicators and applied a correlation threshold of $0.5$, thus considering both strong synergies (correlations $> 0.5$) and strong trade-offs (correlations $< -0.5$). Their results showed that strong negative correlations were absent in the four global regions (OECD, East Asia, Latin America, and MENA). Only one was observed in Sub-Saharan Africa, and two in South Asia. Since the net impact of trade-offs appeared to be weak, they used centrality measures to identify the most influential indicators in each regional network. However, when we applied this methodological approach to the dataset for India, we observed a substantial number of trade-offs, many with significant negative correlations (below $-0.9$). This limits the conclusions drawn in their study, where such strong trade-offs were rarely observed. This observation suggests that while their method is useful in broader regional contexts, it may not be directly applicable to national decision-making.
    \item Song and Jang~\citep{song2023unpacking} suggested another framework to rank SDG targets, by constructing a network based on the similarity between target keywords using semantic analysis. Their approach only establishes a link when two targets share moderate semantic similarities. However, this method neglects the possibility of trade-offs between targets, which are important for understanding the full behaviour of interactions.
\end{itemize}

The drawbacks observed highlight the need for a more flexible method of identifying
important indicators. The primary contribution of this paper is the development of a data-driven network-theoretic framework to classify SDG indicators into
synergy-dominated and trade-off-dominated categories. By introducing the positive and
negative strength measures, our methodology overcomes key limitations of existing
approaches by not relying on fixed correlation thresholds and, crucially, by explicitly
accounting for both synergies and trade-offs rather than neglecting negative
interactions. This ensures that adverse relationships are not discarded or analysed in
isolation, but are systematically incorporated into the assessment of each indicator.
Consequently, the proposed framework provides a more balanced, realistic, and
country-specific classification of SDG indicators.

In addition to the classification framework, we employ an explanatory logistic regression model to
investigate the structural mechanisms underlying synergy dominance. Rather than
serving as a predictive tool, the model is used to provide interpretative insights
into why certain indicators emerge as synergistic. Using this model, we show that indicators classified as synergy-dominated are typically characterised by a combination of significant number of direct positive interactions and a structurally embedded position that allows synergies to extend beyond immediate neighbours. From a policy perspective, this implies that actions targeting such indicators are more likely to align with broader system-wide progress, as their effects can reinforce multiple indicators both directly and through indirect network pathways.

To demonstrate the practical utility of the proposed framework, we conduct case studies
using SDG indicator data for India and Italy. For each country, we classify indicators
into synergy-dominated and trade-off-dominated categories, identifying those that
generate strong positive spillovers and those that induce widespread negative effects.
The comparative analysis shows that both countries exhibit a predominance of
synergy-dominated indicators, reflecting the presence of mutually reinforcing policy
linkages. However, the indicators characterised by synergy and trade-off dominance vary
across national contexts.

In the case of India, access to piped water (SDG 11 - Sustainable Cities and Communities) is found to be a trade-off-dominated indicator, as investments in this area show negative correlations with several health-related indicators. In addition, indicators under Responsible Consumption and Production (SDG 12) and Climate Action (SDG 13) are also found to be trade-off dominated, highlighting the difficulty of achieving environmental objectives without affecting other goals. In contrast, for Italy, indicators under SDG 12 and SDG 13 are synergy-dominated, reflecting how environmental policies have been integrated into broader development goals. Finally, we compare our findings with sector-specific studies and national reports, which support the consistency of our results and demonstrate the ability of the proposed framework to capture systemic interlinkages across different national contexts.

The structure of this paper is as follows. Section~\ref{3 section2} provides the data source for the SDG indicator scores and describes their structure. Section~\ref{3 section3} introduces the proposed network-based framework, including the definitions of positive and negative strength measures used to classify indicators into synergy-dominated and trade-off-dominated categories, and presents the logistic regression framework employed as an explanatory tool to investigate the structural drivers of synergy dominance. Section~\ref{3 section 4} presents the empirical results, discusses the regression findings, and presents the applications to India and Italy to demonstrate the practical relevance of the proposed approach. Section~\ref{3 section 6} concludes the study by summarising the key findings and discussing its limitations.

\section{Data Description}\label{3 section2}

In our study, we used the dataset from the Sustainable Development Report 2025, which we obtained from the official SDG Index website (\url{https://dashboards.sdgindex.org/downloads}). The dataset includes scores for SDG indicators related to the 17 Sustainable Development Goals across different countries, covering the period from 2000 to 2024. Each indicator value ranges from 0 to 100. A score of 100 means the indicator has been achieved, while a score of 0 indicates the poorest performance observed among all countries. These scores enable a comparison of the SDG indicators.

\section{Proposed Methodology}\label{3 section3}
This section introduces the methodological framework used to classify Sustainable Development Goal (SDG) indicators into synergy-dominated and trade-off-dominated categories. The objective is to identify indicators whose improvement is more likely to generate net positive spillovers across the development system, as opposed to those that create competing pressures. The framework consists of three main components:  
\begin{enumerate}
    \item construction of country-specific indicator networks,  
    \item definition of indicator-level measures capturing the balance between synergies and trade-offs, and
    \item an explanatory logistic regression model used to examine the structural drivers of synergy dominance.
\end{enumerate}

\subsection{Network Construction}

Sustainable development indicators are intrinsically interconnected, as progress in one
domain often influences outcomes in others. To systematically represent these interdependencies,
we model the indicator system of each country as a weighted network.

For each country $k$, we construct a complete weighted graph
\[
G_k = (V_k, E_k, W_k),
\]
where $V_k$ is the set of SDG indicators, $E_k$ contains all unordered pairs
$\{i,j\}$ with $i \neq j$, and $W_k = \{w^{(k)}_{ij}\}$ represents edge weights.

The edge weight between indicators $i$ and $j$ is defined as the Spearman rank correlation
coefficient
\[
w^{(k)}_{ij} = \rho^{(k)}_{ij},
\]
computed over the period 2000-2024. The Spearman rank correlation is chosen because it measures the degree of monotonic relationship between two indicators based on their ranked values. The coefficient value ranges from $-1$ to $+1$. A positive correlation indicates that as one indicator increases, the other tends to increase as well. Such relationships are interpreted as synergies, where improvement in one indicator supports or reinforces progress in the other. Conversely, a negative correlation suggests that as one indicator increases, the other tends to decrease. This represents a trade-off, where progress in one may hinder the progress of the other. When the correlation is equal to zero, it implies that the two indicators do not exhibit a monotonic relationship, meaning that the changes in one do not have an association with the changes in the other.

To make the correlation coefficient reliable, indicators with missing data or constant values during the analysis period (2000–2024) are excluded before the construction of the network.

\subsection{Synergy and Trade-off Strength of an Indicator}
Recognising that each indicator is embedded in a complex web of reinforcing and conflicting relationships, we now introduce a measure that captures their combined systemic effect.

Let $G_{ki}$ denote the star subgraph of node $i$ in country $k$, consisting of node $i$ and
all edges incident to it. This subgraph is then divided into:
\begin{itemize}
    \item $G^+_{ki}$, containing only positive-weighted edges (synergies),
    \item $G^-_{ki}$, containing only negative-weighted edges (trade-offs).
\end{itemize}

\begin{figure}[h!]
    \centering
    \includegraphics[width=1\textwidth]{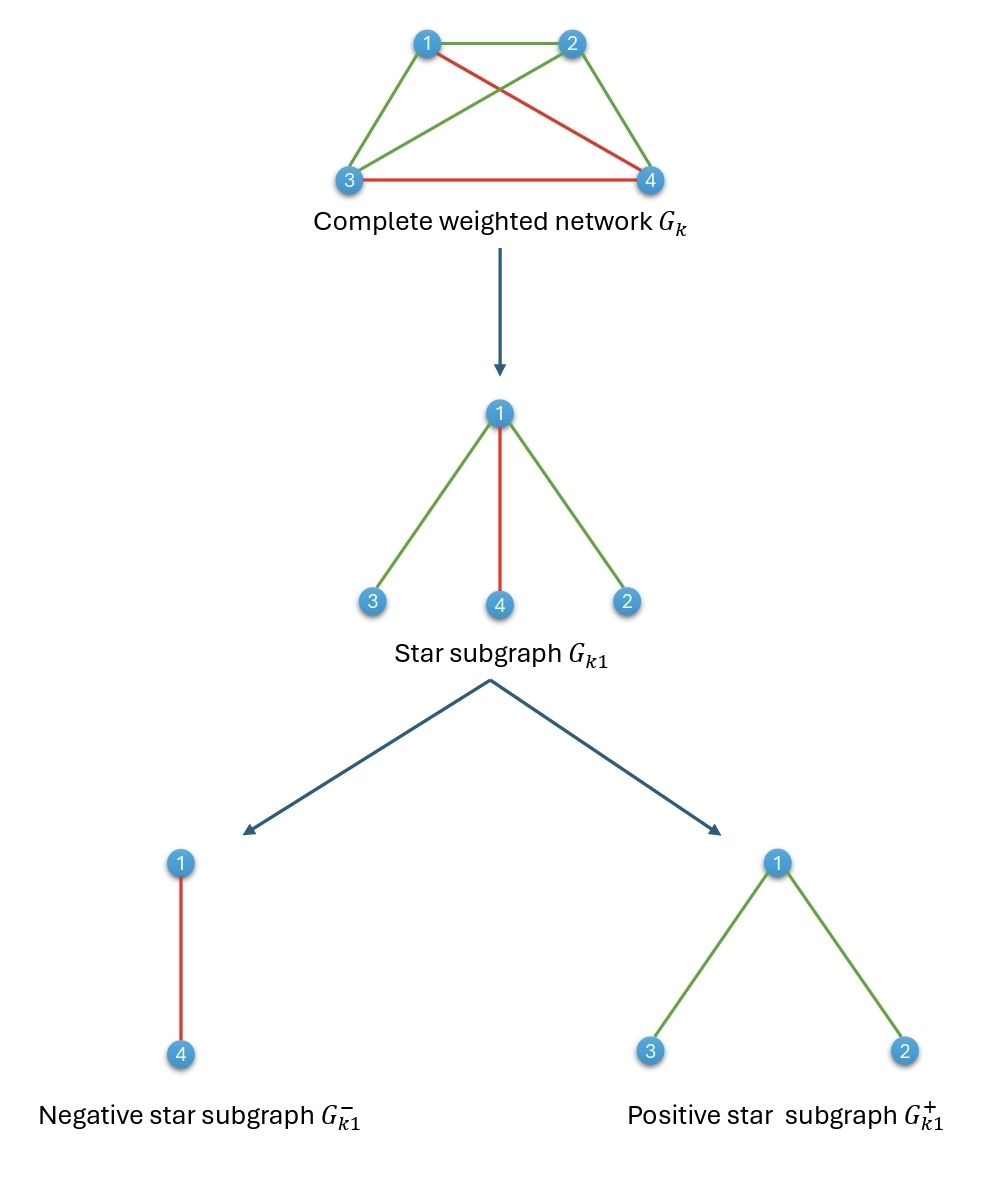}
    \caption{Illustration of a complete weighted network $G_k$ with four nodes, its star subgraph $G_{k1}$ corresponding to node $1$, and its partition into $G^{+}_{k1}$ and $G^{-}_{k1}$ based on postive and negative weights. Green-coloured edges represent the positive weighted edges, while red-coloured edges represent the negative weights.}
    \label{fig:network-illustration}
\end{figure}

Hence, $G_{ki}$ represents the local neighbourhood of node $i$, while $G^+_{ki}$ and $G^-_{ki}$ separately capture its connections corresponding to positive and negative weights. Figure~\ref{fig:network-illustration} illustrates an example of a complete weighted network $G_k$ and the corresponding star subgraph of a chosen node, showing its separation into two weight-based subgraphs.

To evaluate the relative strength of the two classes of connections associated with a node, we now define two measures.

The strength of synergy for an indicator $i$, denoted by $S^+_{ki}$, is defined as:
\begin{equation}\label{3eq:positive-strength}
    S_{ki}^{+} = 
\frac{\sum_{\{i, j\} \in E(G_{ki}^{+})} \rho^{(k)}_{ij}}
     {\sum_{\{i, j\} \in E(G_{ki})} |\rho^{(k)}_{ij}|}
\end{equation}

Similarly, the strength of trade-off for an indicator $i$, denoted by $S^-_{ki}$, is defined as:
\begin{equation}\label{3eq:negative-strength}
 S_{ki}^{-} =  \frac{\sum_{\{i, j\} \in E(G_{ki}^{-})} |\rho^{(k)}_{ij}|}
     {\sum_{\{i, j\} \in E(G_{ki})} |\rho^{(k)}_{ij}|}
\end{equation}

By construction,
\[
0 \le S_{ki}^+ \le 1,\qquad 0 \le S_{ki}^- \le 1,\qquad\text{and}\qquad S_{ki}^+ + S_{ki}^- = 1,
\]
provided $\sum_{\{i, j\} \in E(G_{ki})} |\rho^{(k)}_{ij}| > 0$. 
If the denominator is zero, then $S_{ki}^+$ and $S_{ki}^-$ are undefined, and such indicators can be excluded from the analysis or assigned a zero value by convention. An indicator is said to be synergy dominated if $S^+_{ki} \geq S^-_{ki}$ (or equivalently, if $S^+_{ki} \geq 0.5$), and trade-off dominated otherwise. 

\paragraph{Interpretation and methodological contribution.}
Existing studies on SDG interlinkages typically focus exclusively on synergies, implicitly overlooking the fact that positive and negative interactions often coexist around the same indicator and jointly shape its policy relevance. In contrast, the proposed approach explicitly partitions edge weights into two conceptually meaningful classes. This enables an interpretable indicator-level assessment of the dominance between competing interaction types, which is particularly important in the SDG context where synergies and trade-offs coexist and jointly determine development outcomes. To the best of our knowledge, no existing SDG-focused quantitative framework incorporates such an explicit partitioning to assess the relative dominance of positive and negative interactions at the indicator level, thereby establishing the methodological contribution of the proposed framework.

\subsection{Uncovering the Drivers of Synergy Dominance}
To move beyond a purely rule-based classification and to understand the structural drivers
of synergy dominance, we employ a probabilistic modelling framework. While an indicator
is classified as synergy-dominated when $S^+_{ki} \geq 0.5$, this alone does not
explain why certain indicators exhibit synergistic behaviour. Therefore, we use
logistic regression as an explanatory tool to quantify how specific network properties
influence the likelihood that an indicator is synergy-dominated.

Specifically, we model the conditional probability that an indicator $i$ in country $k$
exhibits synergy dominance as a function of two network-based predictors capturing direct
and indirect effects:

\begin{equation}
\Pr(Y_{ki} = 1 \mid X^d_{ki}, X^h_{ki}) =
\frac{1}{1 + \exp\!\left(-(\beta_0 + \beta_1 X^d_{ki} + \beta_2 X^h_{ki})\right)},
\end{equation}

where:
\begin{itemize}
    \item $Y_{ki}=1$ if indicator $i$ is synergy-dominated and $0$ otherwise,
    \item $X^d_{ki}$ denotes the direct effect, defined as the normalised number of positively weighted edges incident on that node; in particular, $X_{ki}^d$ is calculated by dividing the number of incident edges with positive weights by $n_k - 1$, where $n_k$ denotes the total number of SDG indicators for country $k$. This metric captures the immediate number of positive influences that node $i$ has on the SDG indicator interaction network $G_k$. Normalisation ensures that $X_{ki}^d$ lies between $0$ and $1$, allowing a meaningful comparison between countries with different numbers of indicators.
    \item $X^h_{ki}$ denotes the indirect effect, measured using harmonic
    centrality computed on the subnetwork $G_k^{\text{strong}}$, where $G_k^{\text{strong}}$ is an unweighted subgraph 
     of $G_k$, formed by retaining all vertices of $G_k$ but including only those edges whose weights are at least $0.8$. The harmonic centrality of a node $v_i \in V_k$ quantifies its accessibility from other nodes in the network and is defined as
\[
C_H(v_i) = \frac{1}{n_k-1} \sum_{v_i \ne v_j} \frac{1}{d(v_i, v_j)},
\]
where $d(v_i, v_j)$ denotes the distance\footnote{The path between two nodes $v_i$ and $v_j$ is a sequence of distinct vertices $v_i=v_0,v_1,\dots,v_p=v_j$ such that each consecutive vertices $\{v_{q-1}, v_q\}$ is an edge in $E$. Here $p$ is the length of the path. The distance $d(v_i, v_j)$ between two nodes $v_i$ and $v_j$ is the length of the shortest such path. If there is no path from $v_i$ to $v_j$, then the distance between them is defined to be infinite.} between $v_i$ and $v_j$. Hence, a higher value of $X_{ki}^h$ indicates that node $i$ is more easily reachable from other nodes within the strong-interaction network $G_k^{\text{strong}}$.
 The threshold of $0.8$ is chosen to focus specifically on the flow of strong synergies through the network. By keeping only strong edges, we ensure that any indirect synergetic influence captured in the network arises from strong interactions, rather than being affected by weaker synergies. Since the harmonic centrality of a node $i$ calculates the sum of reciprocals of shortest-path distances from every other node to node $i$, nodes that can reach $i$ in fewer steps of strong synergy contribute more to its score. This means that a node that is easily reachable from all other nodes through short sequences of strong connections will have a higher $X_{ki}^h$. 
    \item $\beta_0, \beta_1, \beta_2$ are parameters to be estimated.
\end{itemize}

The objective of this model is not primarily predictive, but
interpretative. The estimated coefficients allow us to assess:

\begin{enumerate}
    \item whether direct positive interactions ($X^d_{ki}$) significantly increase
    the probability of synergy dominance;
    \item whether the broader structural position of an indicator in the network
    ($X^h_{ki}$) provides additional explanatory power beyond direct connections;
    \item the relative importance of local versus global network effects.
\end{enumerate}

Model parameters are estimated using maximum likelihood. The fitted model yields predicted
probabilities:$$\widehat{\Pr}(Y_{ki}= 1 \mid X_{ki}^d, X_{ki}^h) =
\frac{1}{1 + \exp\!\left(-(\widehat\beta_0 + \widehat\beta_1 X^d_{ki} + \widehat\beta_2 X^h_{ki})\right)},$$ where $\widehat{\beta}_0, \widehat{\beta}_1 \text{ and } \widehat{\beta}_2$ are the estimates of $\beta_0, \beta_1 \text{ and } \beta_2,$ respectively.

An indicator is classified as synergy-dominated when the predicted probability exceeds
$0.5$. However, this classification is secondary; the primary contribution of the
model lies in interpreting the estimated coefficients, which reveal how network structure
shapes synergistic behaviour. By adopting this framework, we can formally test the hypothesis that indicators with stronger direct synergies and more central positions within the strong-interaction
network are more likely to be synergy-dominated. This transforms the analysis from a
descriptive classification exercise into a statistically grounded investigation of the
mechanisms governing SDG indicator interdependencies.

\paragraph{Clarification.}
Note that indicator classification relies exclusively on the positive strength measure $S^+_{ki}$. The regression results are used to interpret how network structure shapes the likelihood of synergy dominance.

\section{Results and Discussions}\label{3 section 4}
This section presents the empirical results obtained from the proposed framework. We first examine the structural drivers underlying synergy dominance by analysing how direct and indirect network effects shape the classification of indicators. We then apply the classification methodology to two national case studies, namely India and Italy, to illustrate how synergy-dominated and trade-off-dominated indicators are distributed across the Sustainable Development Goals in different development contexts.

\subsection{Structural Drivers of Synergy}

Following the methodology described in Section~\ref{3 section2}, we constructed a complete weighted
network $G_k$ for each country $k$ included in the study after removing indicators with
constant or missing values over the period 2000--2024. Using
Equation~\ref{3eq:positive-strength}, we computed the positive strength measure $S^+_{ki}$
for every indicator $i$ in each country $k$. Indicators were classified as synergy-dominated if $S^+_{ki} \ge 0.5$, and trade-off-dominated if $S^+_{ki} < 0.5$.

To understand the structural mechanisms underlying this classification, we
subsequently employed the logistic regression model introduced in
Section~\ref{3 section2}. Importantly, the regression model is used here as an explanatory
tool rather than as a classifier. It allows us to quantify how direct and indirect
network effects shape the likelihood of synergy dominance.

Countries were first grouped according to their SDG Index scores~\citep{sachs2025sdr},
which measure overall progress towards achieving the Sustainable Development Goals.
Countries with scores between $80$ and $100$ were classified as best-performing, those between $50$ and $80$ as moderate-performing, and those below $50$ as worst-performing. The full classification is provided in Table~\ref{3 tab:sdg-categories} (see Appendix). From each performance group, $80\%$ of the indicator observations were randomly selected
to form the training dataset for the regression model, while the remaining $20\%$ were
retained for validation. This stratified sampling ensures that the explanatory model is
calibrated on a balanced representation of development contexts.

The normalised direct effect ($X^d_{ki}$) and indirect effect ($X^h_{ki}$) were used as
predictors, while the binary response variable $Y_{ki}$ indicated synergy dominance
based on $S^+_{ki}$. Variance Inflation Factors (VIF) for both predictors were
approximately $1.55$, indicating the absence of multicollinearity
\citep{james2013introduction}.

The estimated logistic regression model is:

\begin{equation}
\widehat{\Pr}(Y_{ki}= 1 \mid X_{ki}^d, X_{ki}^h) =
\frac{1}{1 + \exp(-(-19.2031 + 39.0684 X_{ki}^d + 2.1742 X_{ki}^h))}.
\label{3 eq:logit_model}
\end{equation}

Rather than serving as a prediction tool, this model enables interpretation of how
network structure drives synergy dominance. When both $X^d_{ki}$ and $X^h_{ki}$ are zero,
the predicted probability approaches zero, indicating that indicators with neither
direct nor indirect influence are extremely unlikely to be synergy-dominated. This is
consistent with intuition, as isolated indicators are unlikely to generate systemic
benefits. Since both coefficients are positive, this implies that increases in either direct or indirect effects increase the likelihood of synergy dominance. Holding $X^d_{ki}$ constant, higher
values of $X^h_{ki}$ lead to higher predicted probabilities, demonstrating that indirect
network embeddedness plays a non-redundant role.

\paragraph{Inference on model parameters.}
The coefficient of $X^d_{ki}$ is estimated as $\widehat{\beta}_1 = 39.07$ with a standard error
of $1.47$ and a $p$-value less than $0.001$, providing strong statistical evidence that
direct synergies play a dominant role in determining synergy dominance. The magnitude of
this coefficient indicates that even small increases in the proportion of positive direct
connections substantially raise the probability that an indicator is classified as
synergy-dominated. Moreover, the associated $95\%$ confidence interval
$[36.20,\;41.94]$ is narrow, suggesting a high degree of precision in the estimate. The
extremely small $p$-value strongly rejects the null hypothesis of no effect, confirming
the robustness of this relationship.

In contrast, the coefficient of the indirect effect is estimated as
$\widehat{\beta}_2 = 2.17$ with a standard error of $0.93$ and a $p$-value of $0.019$. Although
this coefficient is considerably smaller in magnitude, it remains statistically
significant at the $5\%$ level, indicating that indicator position also contributes
positively to synergy dominance. However, the larger standard error and the wider $95\%$
confidence interval $[0.36,\;3.99]$ reflect greater uncertainty around the effect size of
$X^h_{ki}$. This suggests that while indirect influence is not negligible, its impact is
more variable across indicators compared to the direct effect.

Overall, these results demonstrate that direct positive interactions constitute the
primary structural driver of synergy dominance, while indirect network embeddedness plays
a secondary but statistically meaningful role. A complete summary of the regression
estimates is reported in Table~\ref{tab:logit_results}.

\begin{table}[h]
\centering
\caption{Logistic regression coefficients with standard errors, 95\% confidence intervals, and $p$-values.}
\label{tab:logit_results}
\begin{tabular}{lcccc}
\hline
\textbf{Predictor} & \textbf{$\widehat{\beta}$} & \textbf{SE} & \textbf{95\% CI} & \textbf{$p$-value} \\
\hline
Intercept   & $-19.20$ & $0.71$ & $[-20.60,\;-17.81]$ & $<0.001$ \\
$X_{ki}^d$  & $39.07$  & $1.47$ & $[36.20,\;41.94]$   & $<0.001$ \\
$X_{ki}^h$  & $2.17$   & $0.93$ & $[0.36,\;3.99]$     & $0.019$ \\
\hline
\end{tabular}
\end{table}

\paragraph{Validation performance.}

\begin{figure}[h]
    \centering
    \includegraphics[width=0.7\textwidth]{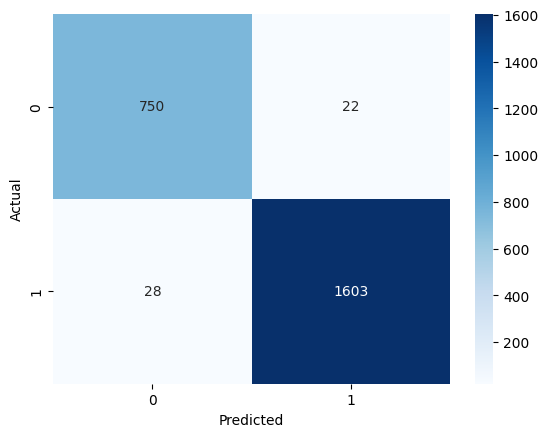}
    \caption{Confusion matrix for the testing dataset}
    \label{3 fig2}
\end{figure}

Using the testing dataset, we assessed model consistency. The model achieved a classification accuracy of $97.92\%$, indicating a high consistency between the structural predictors and the rule-based classification. As shown in Figure~\ref{3 fig2}, which summarises the model's performance on the testing dataset, among all synergy-dominated indicators, $1603$ were correctly classified, while only $28$ were misclassified as trade-off-dominated. Similarly, among trade-off-dominated indicators, $750$ were correctly classified, and $22$ were misclassified as synergy-dominated. This level of consistency supports the use of the identified structural patterns to interpret the systemic role of indicators.

\paragraph{Policy relevance of structural drivers}
An indicator classified as synergy-dominated based on $S^+_{ki}$ is characterised by positive interactions that outweigh its trade-offs, reflecting a net positive systemic influence rather than the total number of connections. The results indicate that indicators with higher values of the direct effect $(X^d_{ki})$ and the indirect effect $(X^h_{ki})$ are more likely to be classified as synergy-dominated. Given the strong consistency observed in the classification results, indicators classified as synergy-dominated based on $S^+_{ki}$ can be viewed as informative entry points for policy design. Since interventions targeting these indicators are more likely to support progress across multiple indicators and allow positive effects to propagate through existing pathways of synergy.

\subsection{Classification of SDG Indicators in India}

\begin{figure}[h]
    \centering
    \includegraphics[width=1\textwidth]{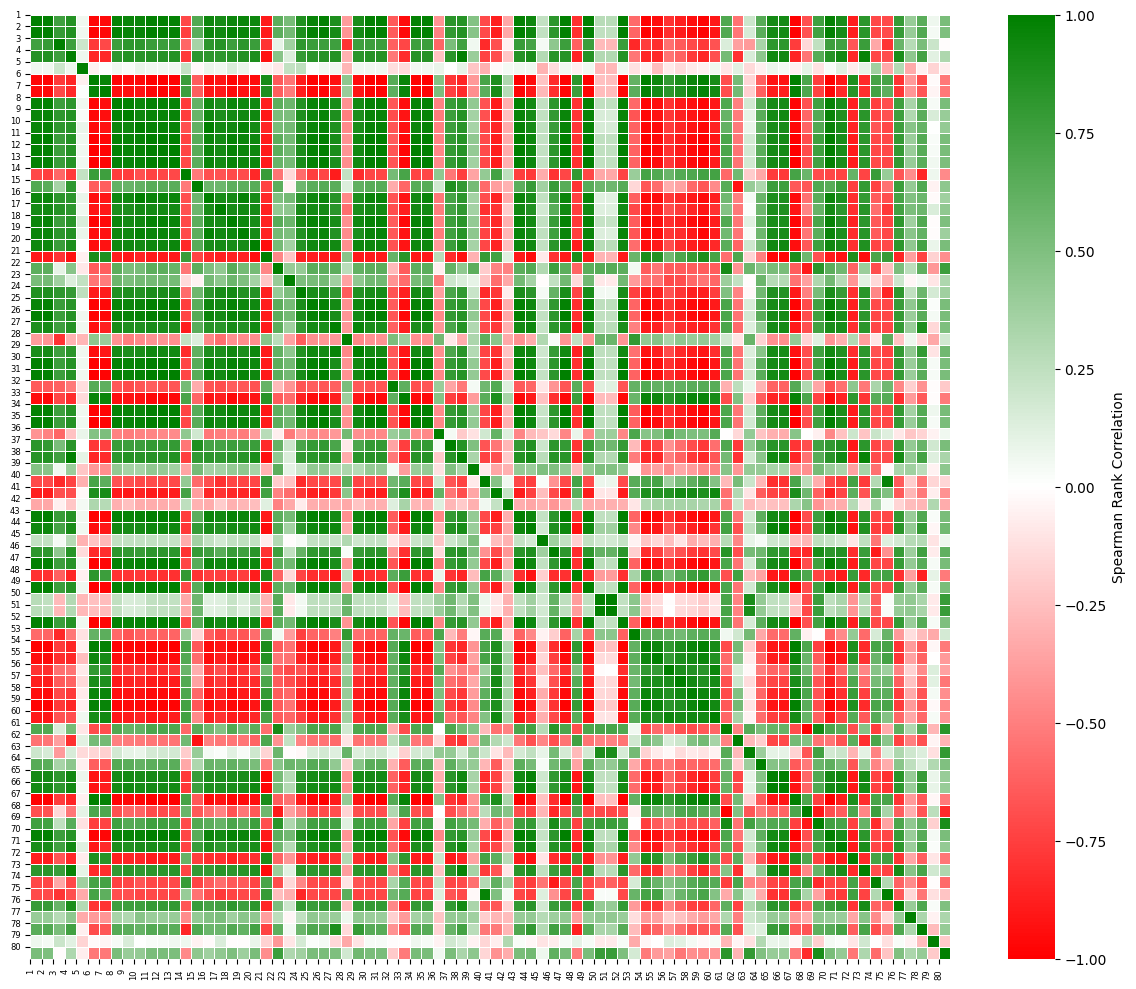}
    \caption{Spearman correlation heatmap of 80 sustainable development indicators for
    India, where dark red indicates strong negative correlation, dark green indicates
    strong positive correlation, and white indicates no correlation.}
    \label{3 correlation heat map}
\end{figure}

To illustrate the country-specific application of the proposed indicator-level classification, we first consider the case of India. All of the nation's SDG indicators were included in the original
dataset. To ensure meaningful Spearman correlations between indicator pairs, we removed
indicators that remained constant throughout the period from 2000 to 2024, as well as
those with missing values during this time. After the data cleaning process, 80
indicators were retained for further analysis. The complete list of the 80 SDG indicators
considered in the study is presented in Table~\ref{3 tab:target_indicators} (see
Appendix). The Spearman rank correlations between the indicators are shown in
Figure~\ref{3 correlation heat map}.

\begin{figure}[h]
    \centering
    \includegraphics[width=0.9\textwidth]{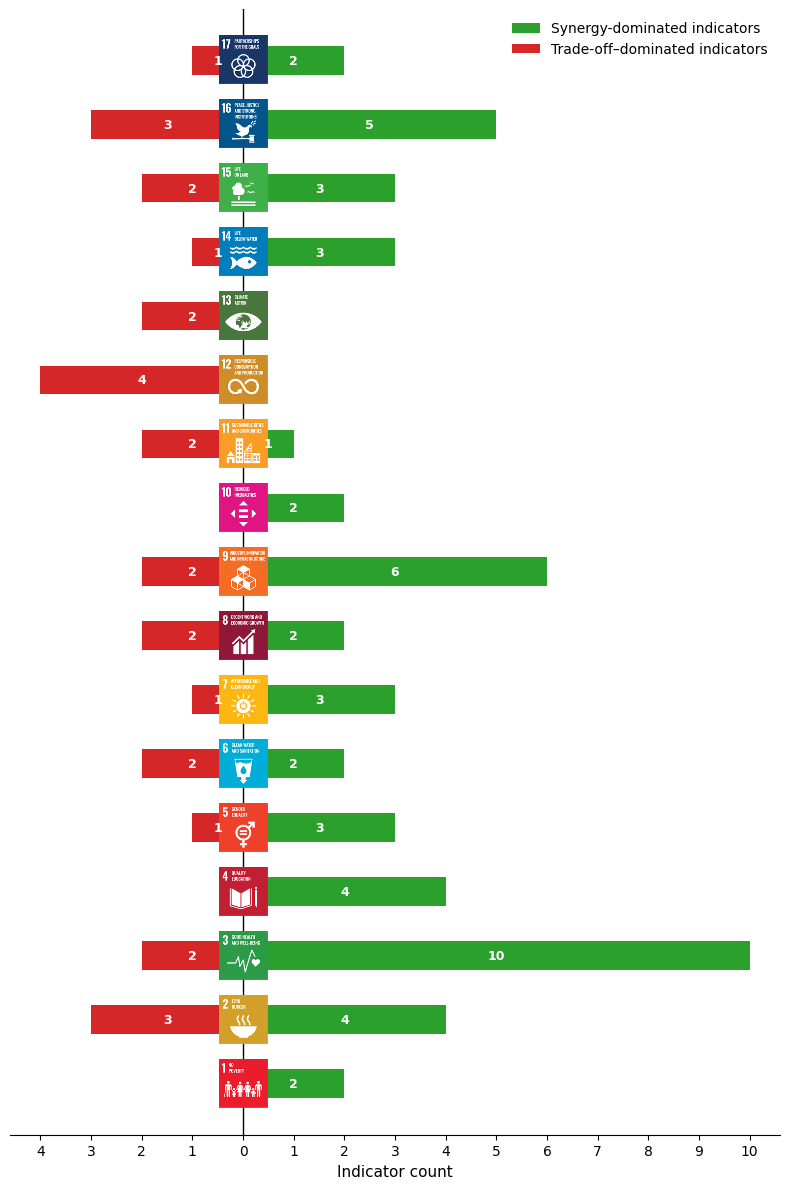}
    \caption{Distribution of synergy and trade-off–dominated SDG indicators across the
    17 Sustainable Development Goals for India. Green bars represent the number of
    indicators that are synergy-dominated, while red bars represent trade-off–dominated
    indicators. Bar lengths indicate indicator counts, shown symmetrically around the
    vertical zero line, with SDG icons aligned at the center for visual reference.}
    \label{3 grouped bar chart}
\end{figure}

Based on the positive strength measure $S^+_{ki}$, $52$ indicators were classified as
synergy-dominated and $28$ as trade-off-dominated. This suggests that the majority of
indicators in the Indian context exhibit synergistic behaviour, showcasing potential for
integrated and mutually reinforcing policy interventions. However, the presence of a
notable number of trade-off-dominated indicators also cautions for careful policy
consideration to avoid unintended negative consequences. To illustrate the distribution
of synergy- and trade-off-dominated indicators across the 17 Sustainable Development
Goals, their counts are presented in Figure~\ref{3 grouped bar chart}.

In Figure~\ref{3 grouped bar chart}, we observe that SDG~3 (Good Health and Well-being)
contains a large number of synergy-dominated indicators, indicating strong positive
interlinkages with other goals and substantial potential for progress across multiple
domains simultaneously. However, SDG~12 (Responsible Consumption and Production) has the
highest number of trade-off-dominated indicators, with none classified as
synergy-dominated. The five indicators of SDG~12 considered in the study include
production-based air pollution, air pollution associated with imports, production-based
nitrogen emissions, and nitrogen emissions associated with imports. This highlights the
challenge in achieving SDG~12, where efforts to reduce pollution and manage waste can
hinder economic or social objectives.

For SDG~13 (Climate Action), both indicators, $CO_2$ emissions from fossil fuel combustion
and cement production, and GHG emissions embodied in imports, exhibit trade-off dominance,
again revealing the tension between economic activity and climate mitigation goals. SDG
11 (Sustainable Cities and Communities) also contains trade-off-dominated indicators. One
such indicator is access to improved water sources through piped connections, which is
classified as trade-off-dominated $(S^+_{ki}=0.314)$. This appears counterintuitive, as
expanding piped water access is generally considered a positive development outcome.
However, Figure~\ref{3 correlation heat map} shows that this indicator exhibits
substantial trade-offs with eleven indicators from SDG~3. Therefore, increasing urban
water infrastructure may improve access, but it can also lead to unintended health
consequences. This observation is supported by an empirical study~\citep{robert2023quality},
which shows that even when households have physical access to piped water, water quality
can be severely degraded due to pollution, poor waste management, or inadequate treatment
infrastructure.

Fishing by trawling or dredging is classified as a synergy-dominated indicator $(S^+_{ki}=0.749)$ for India. This fishing approach damages the seabed by
removing corals, shells, and sediments, which serve as natural shelters and breeding
grounds for many fish species. As a result, fish populations may eventually decline or
fluctuate~\citep{national2002effects}. Reducing such fishing practices therefore supports
ocean ecosystem sustainability. The ocean provides food, jobs, and income for millions,
with some communities entirely dependent on it, particularly in coastal regions. A
healthy and productive ocean contributes to long-term economic growth, public health,
and fish stock maintenance. Coastal communities also benefit from protection against
climate-related hazards such as sea-level rise, erosion, and extreme weather events. In
general, the health of the ocean ecosystem generates wide-ranging social and economic
benefits~\citep{morales2024challenges,pendleton2020we}. 

Overall, the India case study illustrates that SDG indicators differ substantially in their systemic roles, with some supporting broader progress through synergies and others generating notable trade-offs. This underlines the importance of accounting for indicator classifications
when assessing development pathways at the national level.

\subsection{Classification of SDG Indicators in Italy}
\begin{figure}[h]
    \centering
    \includegraphics[width=1\textwidth]{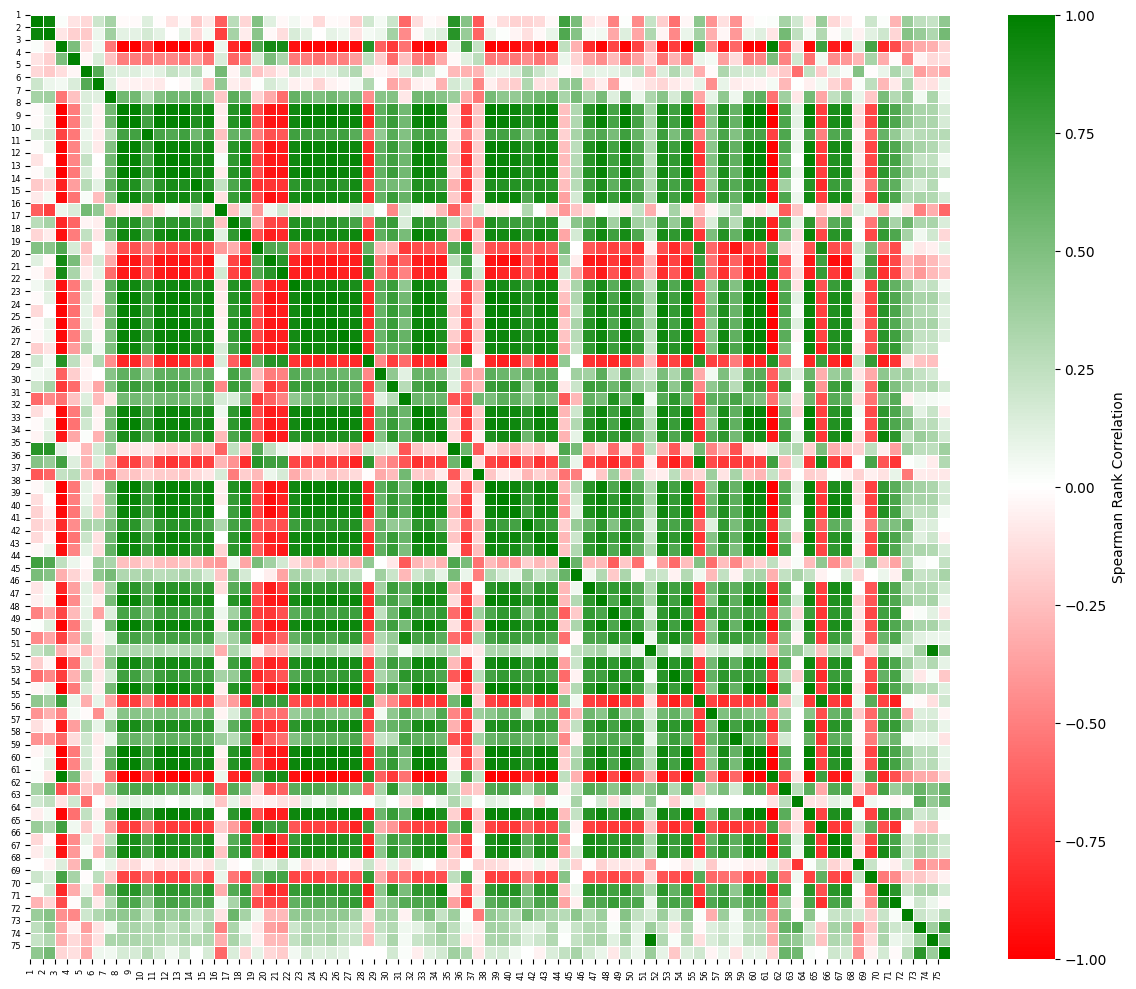}
    \caption{Spearman correlation heatmap of 75 sustainable development indicators for Italy, where dark red indicates strong negative correlation, dark green indicates strong positive correlation, and white indicates no correlation.}
    \label{3 correlation heat map italy}
\end{figure}

To highlight the variations between countries in indicator classifications, we conducted a similar case study for Italy. The initial dataset contained all the SDG indicators available for Italy. To ensure the reliability of the Spearman correlations, indicators that remained constant during the period 2000--2024, or had missing values, were excluded. After the data cleaning process, $75$ indicators were retained for further analysis. The complete list of these $75$ SDG indicators considered in the study, is provided in Table~\ref{3 tab:target_indicators_italy} (see Appendix). The Spearman correlations between indicators are shown in Figure~\ref{3 correlation heat map italy}.

\begin{figure}[h]
    \centering
    \includegraphics[width=0.9\textwidth]{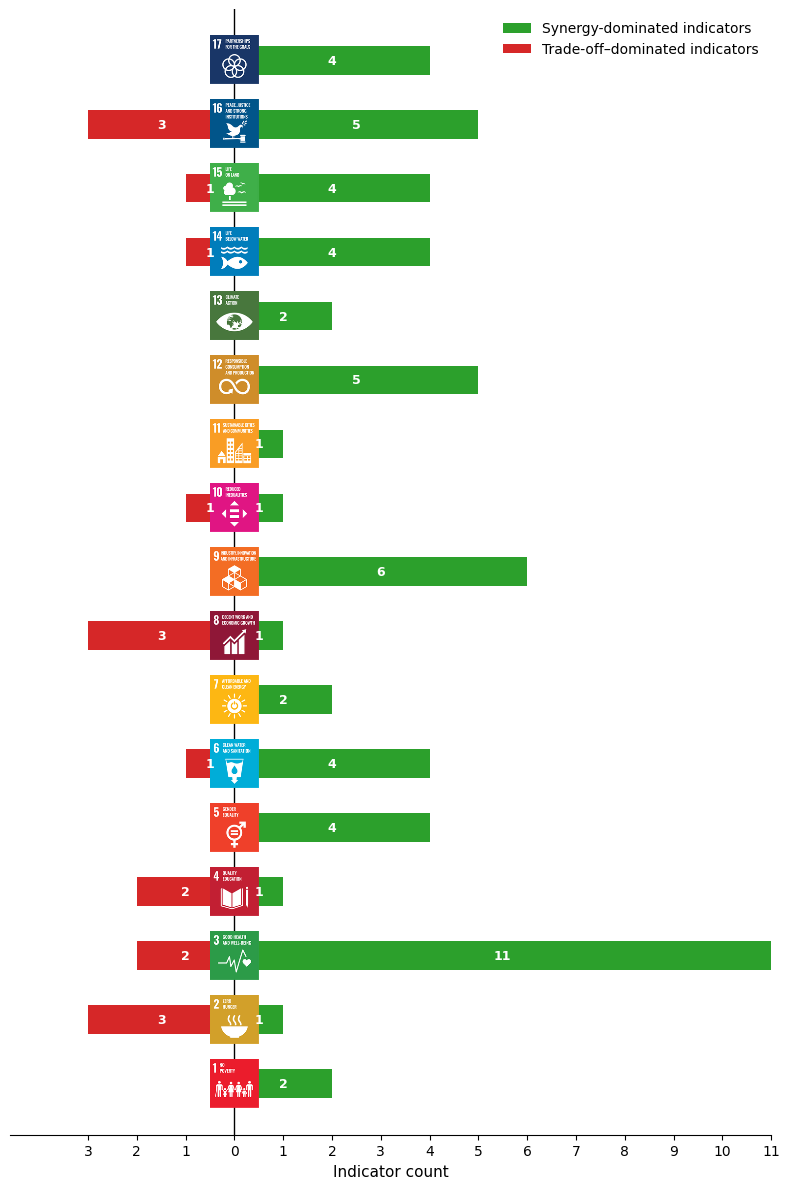}
    \caption{Distribution of synergy and trade-off–dominated SDG indicators across the 17 Sustainable Development Goals for Italy. Green bars represent the number of indicators that are synergy-dominated, while red bars represent trade-off–dominated indicators. Bar lengths indicate indicator counts, shown symmetrically around the vertical zero line, with SDG icons aligned at the center for visual reference.}
    \label{3 grouped bar chart italy}
\end{figure}

Of the $75$ indicators analysed for Italy, $58$ were classified as synergy-dominated and $17$ as trade-off-dominated. This suggests that, similar to India, the majority of Italy’s indicators exhibit synergistic relationships, indicating strong positive interlinkages among development goals. Figure~\ref{3 grouped bar chart italy} presents the distribution of synergy- and trade-off-dominated indicators across the 17 Sustainable Development Goals.

In Figure~\ref{3 grouped bar chart italy}, we note that SDG 3 contains the highest number of synergy-dominated indicators, similar to the case of India. However, unlike India, all five indicators under SDG 12 and both indicators under SDG 13 were classified as synergy-dominated for Italy. These results indicate that improvements in pollution control, nitrogen management, waste reduction, and greenhouse gas mitigation tend to reinforce rather than hinder progress toward other goals. National reports and empirical analyses support this finding. Italy's emissions inventory indicates a steady long-term decrease in the main air pollutants, reflecting the effective implementation of integrated environmental regulation and increased efficiency of industrial processes~\citep{ISPRA2024}. Also, some studies of Italy's waste management performance indicate that achievements under SDG 12 generally have a positive impact on several goals, such as SDG 6 (Clean Water and Sanitation), SDG 7 (Affordable and Clean Energy), and SDG 11 (Sustainable Cities and Communities)~\citep{ram2025achievements}. These developments demonstrate that improvements in the environment support both innovation and public health and contribute to urban sustainability. Thus, they help minimise trade-offs between economic growth and the protection of the environment.

The contrasting patterns observed between India and Italy indicate that the systemic roles of SDG indicators depend strongly on country-specific contexts. Indicators related to environmental regulation and production systems, which appear trade-off-dominated in one context, may become synergy-dominated in another. These findings suggest that indicator-level assessments should be grounded in national data to better reflect local development realities.

\subsection{Global Distribution of Synergy-Dominated Indicators}

\begin{figure}[h]
    \centering
    \includegraphics[width=1\textwidth]{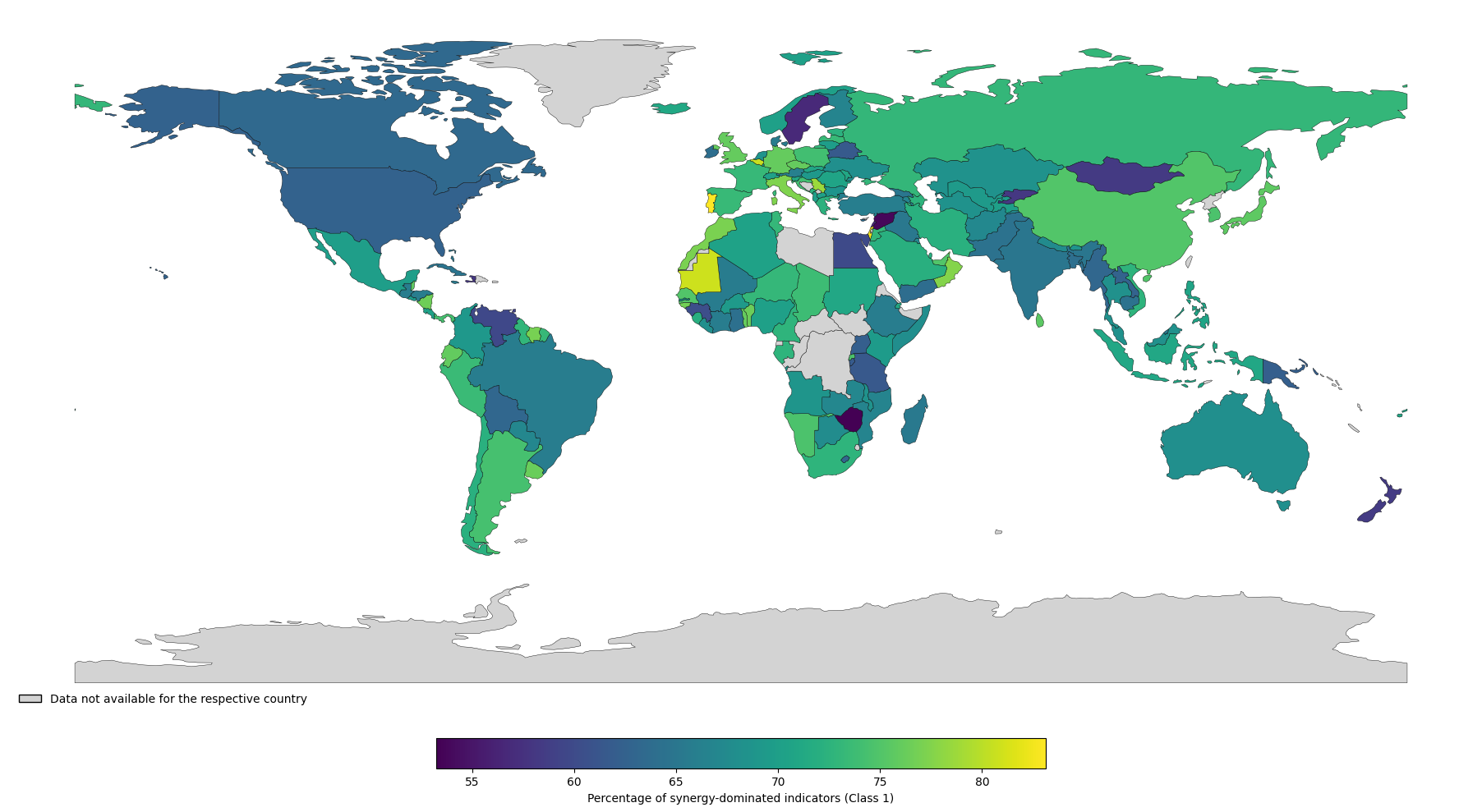}
    \caption{Percentage of synergy-dominated SDG indicators across countries. The colour scale ranges from purple to yellow, where purple shades correspond to countries with lower percentages of synergy-dominated indicators (around $50\%$), blue–green shades indicate intermediate levels, and yellow shades represent countries with high proportions of synergy-dominated indicators (above $80\%$). Grey regions
    indicate countries for which data are not available.}
    \label{3fig:world_synergy_map}
\end{figure}

We now examine the percentage of synergy-dominated indicators across countries around the world. For each country, Figure~\ref{3fig:world_synergy_map} presents the share of indicators classified as synergy-dominated based on the positive strength measure $(S^+_{ki})$.

Across all countries, a substantial proportion of indicators are classified as synergy-dominated, suggesting that net positive interlinkages are a common feature of SDG indicator systems. Even countries in the lowest performance group show a majority of synergy-dominated indicators, indicating that many indicators still exhibit a net positive systemic influence despite weaker overall SDG outcomes. When countries are grouped according to their SDG Index performance, as reported in Table~\ref{3 tab:sdg-categories} (see Appendix), countries in the worst-performing group have, on average, $63.1\%$ synergy-dominated indicators, compared to $68.8\%$ for moderate-performing countries and $72.2\%$ for best-performing countries. This gradual increase suggests that higher-performing countries tend to exhibit more coherent indicator interactions, consistent with more integrated development pathways.

\section{Conclusions}\label{3 section 6}
The Sustainable Development Goals are highly interconnected, and effective implementation requires a careful understanding of the interactions among their underlying indicators. Because policy action is implemented at the indicator level, important insights can be gained by identifying indicators whose improvement generates net positive spillovers, rather than relying solely on goal-level assessments. This study introduces a network-based framework for the quantitative classification of SDG indicators into synergy-dominated and trade-off-dominated categories. The key methodological contribution lies in a flexible, data-driven approach that does not rely on fixed correlation thresholds and explicitly accounts for both positive and negative interactions. By incorporating trade-offs alongside synergies, the framework provides a more balanced and country-specific assessment of indicator-level interdependencies, addressing important gaps in the existing literature.

Beyond classification, the study highlights the structural drivers of synergy dominance. We show that indicators classified as synergy-dominated are typically those that exert direct positive influence on many other indicators, highlighting the importance of local network effects. At the same time, indicators that are more centrally embedded within networks of strong synergies are also more likely to be synergy-dominated, indicating that indirect effects play a complementary role. Together, these findings show that both immediate interactions and the broader network position shape the systemic relevance of SDG indicators, clarifying why certain indicators emerge as strategic points within the development system.

The practical relevance of the framework is demonstrated through national-level applications to India and Italy. In the case of India, SDG~3 contains the largest number of synergy-dominated indicators, while SDG~12 and SDG~13 are characterised primarily by trade-off-dominated indicators. These patterns highlight the challenges of pursuing environmental and production-related objectives without generating competing pressures on health and other development outcomes. In contrast, for Italy, indicators under SDG~12 and SDG~13 are largely synergy-dominated, suggesting that environmental regulation, waste management, and emissions reduction are more effectively integrated with broader development objectives. The comparison between India and Italy underscores that the systemic roles of SDG indicators vary substantially across national contexts, reflecting differences in policy coherence, institutional capacity, and development pathways.

By identifying indicators that generate net positive systemic effects, the proposed framework provides a practical tool for informing policy design and implementation. Indicators classified as synergy-dominated can be interpreted as strategic intervention points, as actions targeting them are more likely to advance multiple goals simultaneously and allow positive effects to propagate through existing synergistic pathways. At the same time, identifying trade-off-dominated indicators highlights areas where careful policy design is required to avoid unintended consequences. In this way, the framework complements existing SDG monitoring tools by translating complex interlinkages into actionable, indicator-level insights.

A key limitation of this study is that it relies on correlation-based relationships, which capture statistical associations between indicators but do not establish causal direction. Since policy decisions are inherently intervention-based, understanding causal mechanisms is crucial for designing effective strategies. Future research should therefore integrate causal inference methods and dynamic modelling approaches to strengthen the identification of leverage points within the SDG system. Extending the framework to incorporate temporal dynamics and causal structure would further enhance its value for evidence-based policymaking and long-term sustainability planning.


\section*{Competing interests}
The authors declare no competing interests.

\section*{Funding}
This research was not supported by any external funding.


\appendix
\section*{Appendix}\label{Appendix}
\begin{table}[ht]
\centering
\renewcommand{\arraystretch}{1.2}
\caption{Categorisation of countries based on SDG Index scores.}
\label{3 tab:sdg-categories}
\begin{tabular}{|p{15cm}|}
\hline
\multicolumn{1}{|c|}{\textbf{Worst-performing countries (SDG Index $<$ 50)}} \\
\hline
Afghanistan, Central African Republic, Chad, Democratic Republic of the Congo, Somalia, South Sudan, Sudan, Yemen\\
\hline
\multicolumn{1}{|c|}{\textbf{Moderate-performing countries (50 $\leq$ SDG Index $<$ 80)}} \\
\hline
Albania, Algeria, Angola, Argentina, Armenia, Australia, Azerbaijan, Bahamas, Bahrain, Bangladesh, Barbados, Belarus, Belize, Benin, Bhutan, Bolivia, Bosnia and Herzegovina, Botswana, Brazil, Brunei Darassalam, Bulgaria, Burkina Faso, Burundi, Cabo Verde, Cambodia, Cameroon, Canada, Chile, China, Colombia, Comoros, Costa Rica, C\^ote d'Ivoire, Cuba, Cyprus, Djibouti, Dominican Republic, Ecuador, Egypt, El Salvador, Eswatini, Ethiopia, Fiji, Gabon, Gambia, Georgia, Ghana, Greece, Guatemala, Guinea, Guinea-Bissau, Guyana, Haiti, Honduras, India, Indonesia, Iran, Iraq, Ireland, Israel, Jamaica, Jordan, Kazakhstan, Kenya, Korea Republic, Kuwait, Kyrgyz Republic, Lao PDR, Lebanon, Lesotho, Liberia, Lithuania, Luxembourg, Madagascar, Malawi, Malaysia, Maldives, Mali, Malta, Mauritania, Mauritius, Mexico, Moldova, Mongolia, Montenegro, Morocco, Mozambique, Myanmar, Namibia, Nepal, New Zealand, Nicaragua, Niger, Nigeria, North Macedonia, Oman, Pakistan, Panama, Papua New Guinea, Paraguay, Peru, Philippines, Qatar, Republic of the Congo, Romania, Russian Federation, Rwanda, S\~{a}o Tom\'{e} and Pr\'{\i}ncipe, Saudi Arabia, Senegal, Serbia, Sierra Leone, Singapore, South Africa, Sri Lanka, Suriname, Switzerland, Syrian Arab Republic, Tajikistan, Tanzania, Thailand, Togo, Trinidad and Tobago, Tunisia, T\"{u}rkiye, Turkmenistan, Uganda, Ukraine, United Arab Emirates, United States, Uruguay, Uzbekistan, Venezuela, Vietnam, Zambia, Zimbabwe\\
\hline
\multicolumn{1}{|c|}{\textbf{Best-performing countries (SDG Index $\geq$ 80)}} \\
\hline
Austria, Belgium, Croatia, Czechia, Denmark, Estonia, Finland, France, Germany, Hungary, Iceland, Italy, Japan, Latvia, Netherlands, Norway, Poland, Portugal, Slovak Republic, Slovenia, Spain, Sweden, United Kingdom\\
\hline
\end{tabular}
\end{table}

\begin{longtable}{clp{11cm}}
\caption{List of SDG Indicators for India}\label{3 tab:target_indicators} \\
\hline
No. & SDG & Indicator Description \\
\hline
\endfirsthead

\hline
No. & SDG & Indicator Description \\
\hline
\endhead

\hline
\endfoot

\hline
\endlastfoot

1 & SDG 1 & Poverty headcount ratio at \$2.15/day (2017 PPP, \%) \\
2 & SDG 1 & Poverty headcount ratio at \$3.65/day (2017 PPP, \%) \\
3 & SDG 2 & Prevalence of undernourishment (\%) \\
4 & SDG 2 & Prevalence of stunting in children under 5 years of age (\%) \\
5 & SDG 2 & Prevalence of wasting in children under 5 years of age (\%) \\
6 & SDG 2 & Prevalence of obesity, BMI $\geq 30$ (\% of adult population) \\
7 & SDG 2 & Human Trophic Level (best 2-3 worst) \\
8 & SDG 2 & Cereal yield (tonnes per hectare of harvested land) \\
9 & SDG 2 & Sustainable Nitrogen Management Index (best 0-1.41 worst)  \\
10 & SDG 3 & Maternal mortality rate (per 100,000 live births) \\
11 & SDG 3 & Neonatal mortality rate (per 1,000 live births) \\
12 & SDG 3 & Mortality rate, under-5 (per 1,000 live births) \\
13 & SDG 3 & Incidence of tuberculosis (per 100,000 population) \\
14 & SDG 3 & Age-standardized death rate due to cardiovascular disease, cancer, diabetes, or chronic respiratory disease in adults aged 30–70 years (\%) \\
15 & SDG 3 & Traffic deaths (per 100,000 population) \\
16 & SDG 3 & Life expectancy at birth (years) \\
17 & SDG 3 & Adolescent fertility rate (births per 1,000 females aged 15 to 19)  \\
18 & SDG 3 & Births attended by skilled health personnel (\%) \\
19 & SDG 3 & Surviving infants who received 2 WHO-recommended vaccines (\%) \\
20 & SDG 3 & Universal health coverage (UHC) index of service coverage (worst 0-100 best) \\
21 & SDG 3 & Subjective well-being (average ladder score, worst 0-10 best) \\
22 & SDG 4 & Participation rate in pre-primary organized learning (\% of children aged 4 to 6) \\
23 & SDG 4 & Net primary enrollment rate (\%) \\
24 & SDG 4 & Lower secondary completion rate (\%) \\
25 & SDG 4 & Literacy rate (\% of population aged 15 to 24) \\
26 & SDG 5 & Demand for family planning satisfied by modern methods (\% of females aged 15 to 49) \\
27 & SDG 5 & Ratio of female-to-male mean years of education received (\%) \\
28 & SDG 5 & Ratio of female-to-male labor force participation rate (\%) \\
29 & SDG 5 & Seats held by women in national parliament (\%) \\
30 & SDG 6 & Population using at least basic drinking water services (\%) \\
31 & SDG 6 & Population using at least basic sanitation services (\%) \\
32 & SDG 6 & Freshwater withdrawal (\% of available freshwater resources) \\
33 & SDG 6 & Scarce water consumption embodied in imports ($m^3H_2O$ eq/capita) \\
34 & SDG 7 & Population with access to electricity (\%) \\
35 & SDG 7 & Population with access to clean fuels and technology for cooking (\%) \\
36 & SDG 7 & $CO_2$ emissions from fuel combustion per total electricity output $(MtCO_2/TWh)$  \\
37 & SDG 7 & Renewable energy share in total final energy consumption (\%) \\
38 & SDG 8 & Adults with an account at a bank or other financial institution or with a mobile-money-service provider (\% of population aged 15 or over) \\
39 & SDG 8 & Unemployment rate (\% of total labor force, ages 15+) \\
40 & SDG 8 & Fundamental labor rights are effectively guaranteed (worst 0–1 best) \\
41 & SDG 8 & Fatal work-related accidents embodied in imports (per million population) \\
42 & SDG 9 & Rural population with access to all-season roads (\%) \\
43 & SDG 9 & Population using the internet (\%) \\
44 & SDG 9 & Mobile broadband subscriptions (per 100 population) \\
45 & SDG 9 & Logistics Performance Index: Infrastructure Score (worst 1–5 best) \\
46 & SDG 9 & The Times Higher Education Universities Ranking: Average score of top 3 universities (worst 0-100 best) \\
47 & SDG 9 & Articles published in academic journals (per 1,000 population) \\
48 & SDG 9 & Expenditure on research and development (\% of GDP) \\
49 & SDG 9 & Total patent applications by applicant's origin (per million population) \\
50 & SDG 10 & Gini coefficient \\
51 & SDG 10 & Palma ratio \\
52 & SDG 11 & Proportion of urban population living in slums (\%) \\
53 & SDG 11 & Annual mean concentration of PM2.5 ($\mu g / m^3$) \\
54 & SDG 11 & Access to improved water source, piped (\% of urban population) \\
55 & SDG 12 & Production-based air pollution (DALYs per 1,000 population) \\
56 & SDG 12 & Air pollution associated with imports (DALYs per 1,000 population) \\
57 & SDG 12 & Production-based nitrogen emissions (kg/capita) \\
58 & SDG 12 & Nitrogen emissions associated with imports (kg/capita) \\
59 & SDG 13 & $CO_2$ emissions from fossil fuel combustion and cement production (tCO2/capita) \\
60 & SDG 13 & GHG emissions embodied in imports ($tCO_2$/capita) \\
61 & SDG 14 & Ocean Health Index: Clean Waters score (worst 0-100 best) \\
62 & SDG 14 & Fish caught from overexploited or collapsed stocks (\% of total catch) \\
63 & SDG 14 & Fish caught by trawling or dredging (\%) \\
64 & SDG 14 & Fish caught that are then discarded (\%) \\
65 & SDG 15 & Mean area that is protected in terrestrial sites important to biodiversity (\%) \\
66 & SDG 15 & Mean area that is protected in freshwater sites important to biodiversity (\%) \\
67 & SDG 15 & Red List Index of species survival (worst 0-1 best)  \\
68 & SDG 15 & Permanent deforestation (\% of forest area, 3-year average) \\
69 & SDG 15 & Imported deforestation ($m^2$/capita) \\
70 & SDG 16 & Homicides (per 100,000 population) \\
71 & SDG 16 & Crime is effectively controlled \\
72 & SDG 16 & Unsentenced detainees (\% of prison population) \\
73 & SDG 16 & Corruption Perceptions Index (worst 0-100 best) \\
74 & SDG 16 & Press Freedom Index (worst 0-100 best) \\
75 & SDG 16 & Access to and affordability of justice (worst 0–1 best) \\
76 & SDG 16 & Timeliness of administrative proceedings (worst 0 - 1 best) \\
77 & SDG 16 & Expropriations are lawful and adequately compensated (worst 0 - 1 best) \\
78 & SDG 17 & Government spending on health and education (\% of GDP) \\
79 & SDG 17 & Other countries: Government revenue excluding grants (\% of GDP) \\
80 & SDG 17 & Statistical Performance Index (worst 0-100 best) \\

\end{longtable}

\begin{longtable}{clp{11cm}}
\caption{List of SDG Indicators for Italy}\label{3 tab:target_indicators_italy} \\
\hline
No. & SDG & Indicator Description \\
\hline
\endfirsthead

\hline
No. & SDG & Indicator Description \\
\hline
\endhead

\hline
\endfoot

\hline
\endlastfoot

1 & SDG 1 & Poverty headcount ratio at \$2.15/day (\%) \\
2 & SDG 1 & Poverty headcount ratio at \$3.65/day (\%) \\
3 & SDG 2 & Prevalence of obesity, BMI $\geq 30$ (\% of adult population) \\
4 & SDG 2 & Human Trophic Level (best 2–3 worst) \\
5 & SDG 2 & Cereal yield (tonnes per hectare of harvested land) \\
6 & SDG 2 & Sustainable Nitrogen Management Index (best 0–1.41 worst) \\
7 & SDG 3 & Maternal mortality ratio (per 100,000 live births) \\
8 & SDG 3 & Neonatal mortality rate (per 1,000 live births) \\
9 & SDG 3 & Mortality rate, under-5 (per 1,000 live births) \\
10 & SDG 3 & Incidence of tuberculosis (per 100,000 population) \\
11 & SDG 3 & New HIV infections (per 1,000 uninfected population, all ages) \\
12 & SDG 3 & Age-standardized death rate due to cardiovascular disease, cancer, diabetes, or chronic respiratory disease in adults aged 30–70 years (\%) \\
13 & SDG 3 & Traffic deaths (per 100,000 population) \\
14 & SDG 3 & Life expectancy at birth (years) \\
15 & SDG 3 & Adolescent fertility rate (births per 1,000 females aged 15–19) \\
16 & SDG 3 & Births attended by skilled health personnel (\%) \\
17 & SDG 3 & Surviving infants who received 2 WHO-recommended vaccines (\%) \\
18 & SDG 3 & Universal health coverage (UHC) index of service coverage (worst 0–100 best) \\
19 & SDG 3 & Subjective well-being (average ladder score, worst 0–10 best) \\
20 & SDG 4 & Participation rate in pre-primary organized learning (\% of children aged 4 to 6) \\
21 & SDG 4 & Net primary enrollment rate (\%) \\
22 & SDG 4 & Lower secondary completion rate (\%) \\
23 & SDG 5 & Demand for family planning satisfied by modern methods (\% of females aged 15–49) \\
24 & SDG 5 & Ratio of female-to-male mean years of education received (\%) \\
25 & SDG 5 & Ratio of female-to-male labor force participation rate (\%) \\
26 & SDG 5 & Seats held by women in national parliament (\%) \\
27 & SDG 6 & Population using at least basic drinking water services (\%) \\
28 & SDG 6 & Population using at least basic sanitation services (\%) \\
29 & SDG 6 & Freshwater withdrawal (\% of available freshwater resources) \\
30 & SDG 6 & Anthropogenic wastewater that receives treatment (\%) \\
31 & SDG 6 & Scarce water consumption embodied in imports ($m^3H_2O$ eq/capita) \\
32 & SDG 7 & $CO_2$ emissions from fuel combustion per total electricity output $(MtCO_2/TWh)$ \\
33 & SDG 7 & Renewable energy share in total final energy consumption (\%) \\
34 & SDG 8 & Adults with an account at a bank or other financial institution or with a mobile-money-service provider (\% of population aged 15 or over) \\
35 & SDG 8 & Unemployment rate (\% of total labor force, ages 15+) \\
36 & SDG 8 & Fundamental labor rights are effectively guaranteed (worst 0–1 best) \\
37 & SDG 8 & Fatal work-related accidents embodied in imports (per million population) \\
38 & SDG 9 & Population using the internet (\%) \\
39 & SDG 9 & Mobile broadband subscriptions (per 100 population) \\
40 & SDG 9 & Logistics Performance Index: Infrastructure Score (worst 1–5 best) \\
41 & SDG 9 & Articles published in academic journals (per 1,000 population) \\
42 & SDG 9 & Expenditure on research and development (\% of GDP) \\
43 & SDG 9 & Total patent applications by applicant's origin (per million population) \\
44 & SDG 10 & Gini coefficient \\
45 & SDG 10 & Palma ratio \\
46 & SDG 11 & Annual mean concentration of PM2.5 ($\mu g/m^3$) \\
47 & SDG 12 & Production-based air pollution (DALYs per 1,000 population) \\
48 & SDG 12 & Air pollution associated with imports (DALYs per 1,000 population) \\
49 & SDG 12 & Production-based nitrogen emissions (kg/capita) \\
50 & SDG 12 & Nitrogen emissions associated with imports (kg/capita) \\
51 & SDG 12 & Exports of plastic waste (kg/capita) \\
52 & SDG 13 & $CO_2$ emissions from fossil fuel combustion and cement production ($tCO_2$/capita) \\
53 & SDG 13 & GHG emissions embodied in imports ($tCO_2$/capita) \\
54 & SDG 14 & Mean area that is protected in marine sites important to biodiversity (\%) \\
55 & SDG 14 & Ocean Health Index: Clean Waters score (worst 0–100 best) \\
56 & SDG 14 & Fish caught from overexploited or collapsed stocks (\% of total catch) \\
57 & SDG 14 & Fish caught by trawling or dredging (\%) \\
58 & SDG 14 & Fish caught that are then discarded (\%) \\
59 & SDG 15 & Mean area that is protected in terrestrial sites important to biodiversity (\%) \\
60 & SDG 15 & Mean area that is protected in freshwater sites important to biodiversity (\%) \\
61 & SDG 15 & Red List Index of species survival (worst 0–1 best) \\
62 & SDG 15 & Permanent deforestation (\% of forest area, 3-year average) \\
63 & SDG 15 & Imported deforestation ($m^2$/capita) \\
64 & SDG 16 & Homicides (per 100,000 population) \\
65 & SDG 16 & Crime is effectively controlled (worst 0–1 best) \\
66 & SDG 16 & Unsentenced detainees (\% of prison population) \\
67 & SDG 16 & Corruption Perceptions Index (worst 0–100 best) \\
68 & SDG 16 & Press Freedom Index (worst 0–100 best) \\
69 & SDG 16 & Access to and affordability of justice (worst 0–1 best) \\
70 & SDG 16 & Timeliness of administrative proceedings (worst 0–1 best) \\
71 & SDG 16 & Expropriations are lawful and adequately compensated (worst 0–1 best) \\
72 & SDG 17 & Government spending on health and education (\% of GDP) \\
73 & SDG 17 & International concessional public finance, including official development assistance (\% of GNI) \\
74 & SDG 17 & Corporate Tax Haven Score (best 0–100 worst) \\
75 & SDG 17 & Statistical Performance Index (worst 0–100 best) \\

\end{longtable}

\bibliographystyle{apalike}  
\bibliography{references}  

\end{document}